# Breaking the Resolution limit in Photoacoustic Imaging using Positivity and Sparsity

Peter Burgholzer[a]
Johannes Bauer-Marschallinger[a]
Mike Hettich[a]
Markus Haltmeier[b]

[a] Research Center for Non Destructive Testing (RECENDT), Linz, Austria; peter.burgholzer@recendt.at
[b] Department of Mathematics, University of Innsbruck, Austria

**Abstract:** The spatial resolution in photoacoustic imaging is essentially limited by acoustic attenuation, which can be numerically compensated only up to a theoretical limit. The physical background for this "ill-posedness" is the second law of thermodynamics: the loss of information is equal to the entropy production, which is the energy decay of the attenuated wave divided by the temperature. As acoustic attenuation increases with higher frequencies, a cut-off frequency can be determined, where the information content for that frequency gets so low that it cannot be distinguished from equilibrium distribution within a certain statistical significance. This cut-off frequency can be determined also by setting the amplitude of the attenuated signal in frequency domain equal to the noise-level. Compensating for acoustic attenuation requires to solve an ill-posed inverse problem, where an adequate regularization parameter is the cut-off frequency, when the acoustic wave amplitude is damped just below the noise level. If additional information, such as positivity or sparsity is used, this theoretical resolution limit can be overcome. This is experimentally demonstrated for the propagation of planar acoustic waves in fat tissue, which are induced by short laser pulses and measured by piezoelectric transducers. For fatty porcine tissue the frequency dependent acoustic attenuation was measured. This was used to invert the problem and by using additional information, in the form of positivity and sparsity (Douglas-Rachford splitting algorithm) the resolution could be enhanced significantly compared to the limit given by the cut-off frequency from attenuation through 20 mm of porcine fat tissue.

# Introduction

Photoacoustic (or optoacoustic) imaging uses the thermo-elastic expansion following a rapid temperature rise after illumination of light absorbing structures within a semitransparent and turbid material, such as a biological tissue. It allows acoustic resolution with simultaneous optical absorption contrast and enables to detect hemoglobin, lipids, water and other light-absorbing chromophores, with greater penetration depth than with purely optical imaging modalities that rely on ballistic photons [1,2,3]. In photoacoustic tomography, the temporal evolution of the acoustic pressure field is sampled using an array of ultrasound detectors placed on the sample surface or by moving a single detector across the sample surface. From the measured pressure signals, images of the optical absorption within the tissue are reconstructed by solving an inverse source problem [3,4,5].

In this work, the achievable spatial resolution for photoacoustic imaging in porcine fat tissue is investigated. At depths larger than the range of the ballistic photons, i.e. more than a few hundreds of microns in tissue, light is scattered several times and the spatial resolution in photoacoustic imaging is limited by acoustic attenuation, which is caused by acoustic absorption, dispersion, and scattering. The spatial resolution degrades with increasing depth because higher acoustic frequencies, which have smaller wavelengths and allow a better resolution, are stronger attenuated than lower frequencies. As a rule of thumb, the ratio of the imaging depth to the best spatial resolution is roughly constant and has a value of 200 [3]. Although acoustic attenuation defines the ultimate spatial resolution limit, other factors such as detector bandwidth, element size and the area over which the acoustic signals are recorded at the sample surface – the detection aperture – can be limiting factors in practice [2]. These technical limitations can be avoided in principle – or at least can be reduced.

Mathematically, the compensation of frequency-dependent acoustic attenuation is an ill-posed inverse problem, where the cut-off frequency is an adequate regularization parameter. The physical reason for the ill-posedness is the second law of thermodynamics: acoustic attenuation is an irreversible process and the entropy production, which is the energy decay during wave propagation due to attenuation divided by the temperature, is equal to the information loss for the reconstructed image [6]. This information loss due to entropy production cannot be compensated mathematically. As the information content of the reconstructed image strongly correlates with the spatial resolution, this results in a fundamental resolution limit due to thermodynamic principles. Non-equilibrium thermodynamics describes the connection between entropy production and information loss, e.g. [7], and as we have already elaborated for heat diffusion [8,9], we could determine a cut-off frequency also for damped acoustic waves in water [10]. Here, this approach is generalized for a frequency dependent acoustic attenuation described by a power-law using a general exponent which might be different from two in water.

In frequency space, the information content of wave components with frequencies above that cut-off frequency is so low that they cannot be statistically distinguished from the equilibrium distribution. This is equivalent that the acoustic wave amplitude on the sample surface is damped just below the noise level [10]. Consequently, the spatial resolution limit becomes diffraction limited and according to Nyquist it is half of the wavelength at this cut-off frequency [10]. It is not a fortunate coincidence that the same cut-off frequency can be determined from entropy production and from noise-fluctuations, but stems from the fluctuation-dissipation-relation described in statistical physics [6,8,9]. To reach this thermodynamic resolution limit for compensation of acoustic attenuation experimentally, it is necessary to measure the broadband ultrasonic attenuation parameters of tissues or liquids very accurately [11] and to evaluate the existing mathematical models to get an adequate description of attenuation [12]. For the used porcine fat tissue in the measured frequency range, a power law very well describes the dependence of the attenuation on frequency. We emphasize, however, that the proposed method for evaluating the principle resolution limit is applicable to any acoustic attenuation model described by a complex wave number.

In this work it is demonstrated, that even for "ordinary" photoacoustic measurements using one image with homogeneous illumination and no moving particles or droplets this resolution limit can be exceeded significantly by taking non-negativity and sparsity information as an additional knowledge into account. Therefore, this resolution enhancement can be reached for all existing conventional photoacoustic tomography or acoustic resolution microscopy set-ups without additional time due to multiple measurements. This was partially inspired by works from geophysics, where the attenuation of seismic signals was compensated [13,14]. For one-dimensional (1D) photoacoustic pressure pulses, non-negativity is evident, as the initial pressure generated by optical heating is always positive. Propagation in 3D gives negative pressure components, but transforming to spherical projections, which is the time integral of the measured acoustic pressure, yields again positivity for all measured signals [15].

For regularization two different methods are compared: the truncated singular value decomposition (T-SVD) method [16], which allows also negative values and does not use sparsity, and the "Douglas-Rachford splitting algorithm" (DR algorithm), which uses positivity and sparsity [17,18]. By taking positivity and sparsity into account, the DR algorithm can enhance resolution significantly. Experimentally, this is demonstrated in 1D, but as mentioned in section 1.1, the same resolution enhancement is expected in axial direction in 3D by using spherical projections. In lateral direction, due to limited angle effects the resolution enhancement is less, but could be partly compensated by using weight factors [19]. Here, the resolution limit is defined as half of the wavelength of the cut-off frequency, where the acoustic signal at the surface is damped just below the noise level. For linear reconstruction methods, this is equivalent to the width of the reconstruction of a point-source signal, which is the point-spread function. For nonlinear methods, such as the used DR algorithm, the blurring depends on the number of iterations and the acoustic signal itself, and therefore resolution cannot be defined as the width of the point-spread function any more, but how well distinct sources can be distinguished from each other in the reconstructed image. Therefore, resolution has to be defined as a measure of localization and separation.

# Experimental Setup

Ultrasonic waves were generated, attenuated, and detected experimentally in a setup sketched in Figure 1. To excite strong and broadband ultrasonic signals, a laser ultrasound method is employed [20,21]. Short nanosecond laser pulses are directed on a silicon wafer. A cylindrical region of the wafer, ranging from the surface to a depth of a few microns and with a diameter determined by the laser beam, abruptly heats up by optical absorption. The subsequent thermoelastic expansion of the heated volume leads to the emission of ultrasonic waves mainly directed perpendicular to the surface of the wafer. Porcine subcutaneous fat tissue in the propagation path of the sound waves leads to acoustic attenuation due to absorption and scattering. After exiting the tissue, the sound waves propagate through water and are detected by a piezoelectric transducer. The resulting electrical signals are amplified (5073PR-40-E, Olympus NDT Inc., Waltham, MA) and sampled by an oscilloscope (DSO 5043, 300 MHz; Agilent Technologies Inc., Santa Clara, CA).

Two mounting flanges hold the parts by a small axial force applied by two screws to keep a well-defined distance to the piezoelectric transducer of 6 mm and 20 mm. To get simultaneously two signals slightly shifted in time through the fatty tissue, a small step is fabricated by ion milling near the beam center of the silicon wafer. By using an unfocused piezoelectric transducer (V358-SU, Panametrics, Waltham, MA) the influence of possibly occurring local variation of the attenuation is decreased. If the silicon wafer shows a small step near the beam center the two slightly shifted signals in time through the fatty tissue overlap. The transducer has a center frequency of 50.6 MHz, a diameter of the sensing element of 6.35 mm and a -6 dB bandwidth of 81.2 %. Within the observed bandwidth, the acoustic attenuation occurring in the short paths of water in the setup can be safely neglected compared to the strong attenuation caused by the relative thick layers of fat [10,22]. Therefore, different water path lengths do not change the measured signals. The mechanical set-up is designed in a way to ensure parallel alignment of the silicon wafer, the fat tissue and the piezoelectric transducer. Inclination of these components would lead to unwanted signal losses due to refraction and reflection. Also, since the waves show plane wave behavior, an inclination of the wavefront relative to the sensing element of the transducer will lead to a misleading low-pass-filtering of the signals. A detailed description of the setup can be found in [11]. The ultrasound-generating optical pulses have a diameter of 6 mm and are emitted by a frequency-doubled Nd:YAG laser (Continuum Surelite, 20 Hz repetition rate, 6 ns pulse duration, 532 nm center wavelength). The used pulse energies ranged from 12 mJ to 65 mJ.

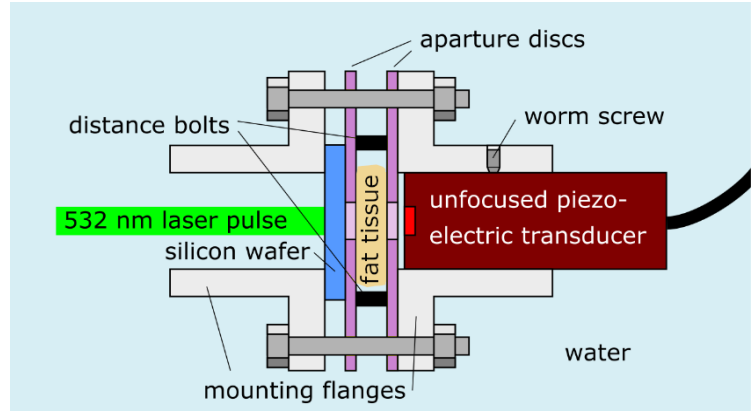

Figure 1. Setup for the generation and detection of acoustic plane waves attenuated by porcine fat tissue. Abrupt local heating of a silicon wafer by short laser pulses leads to the emission of strong broadband ultrasonic plane waves. Porcine subcutaneous fat tissue in the propagation path induces frequency dependent attenuation of the acoustic signals. The fat tissue is fixated between two aperture disks applying a small axial force on the tissue. Distance bolts with 6 mm and respectively 20 mm ensure two precise lengths of the attenuation path. For these two lengths, attenuated acoustic plane waves were detected by an unfocused piezoelectric transducer. Image from https://doi.org/10.3390/jimaging5010013 by Burgholzer et al. was edited and is used under CC BY 4.0.

## Compensation of Acoustic Attenuation and Resolution Limits

As mentioned in section 1, the compensation of acoustic attenuation in photoacoustic imaging can be reduced to a 1D-problem. The relation between the ideal acoustic wave pressure $p_{ideal}(r,t)$ without attenuation and the attenuated wave pressure $p(r,t)$ at the same 1D, 2D, or 3D location $r$ is identical for all dimensions, as shown by Ammari et al. in a compact way [23], which reads in frequency domain:

$$\tilde{p}(r,\omega) = \frac{\omega}{c_0 K(\omega)} \widetilde{p_{ideal}}(r, c_0 K(\omega)), \qquad (1)$$

where $\omega$ is the angular frequency, and the tilde indicates the Fourier transformation of the time signals $p_{ideal}(r,t)$ and $p(r,t)$, respectively. The complex wave number $K(\omega) := \omega/c(\omega) + i\,\alpha(\omega)$ describes the frequency dependent dispersion and attenuation with the phase velocity $c(\omega)$ and the attenuation coefficient $\alpha(\omega)$. For the ideal wave without attenuation, the wavenumber $k(\omega) := \omega/c_0$ is real and $c_0$ is the sound velocity for the ideal wave, which shows no dispersion and frequency dependence. Phase velocity and attenuation coefficient are connected by a Kramers-Kronig relationship to guarantee causality [24].

For an acoustic signal, which propagates through an attenuating sample of a defined thickness, the effect of acoustic attenuation in 3D and its compensation was modeled in the frequency domain and 3D by Dean-Ben et al [22]. Here, a similar derivation will be performed in 1D, as we get plane waves according to our experimental setup described in section 2.1. Because the acoustic attenuation in water compared to fat tissue for the bandwidth covered by our detector can be neglected [10,22], the signal measured in water without any fat tissue is taken as the ideal one. With the Helmholtz equation $(\nabla^2 + K(\omega)^2)\tilde{p}(r,\omega) = \delta(r)$ for a uniform attenuating medium and a point source located at the origin $r = 0$ [22], one gets by using the 1D Greens function $\tilde{p}(r,\omega) = \exp(iK(\omega)|r|)/(2iK(\omega))$ the relation between the attenuated and the ideal signal:

$$\tilde{p}(r,\omega) = \frac{\omega}{c_0 K(\omega)} \exp(i\gamma(\omega)|r|)\,\widetilde{p_{ideal}}(r,\omega), \qquad (2)$$

where $\gamma(\omega)$ is the difference of the complex wavenumber for the attenuated wave and the real wave number for the ideal wave:

$$\gamma(\omega) := K(\omega) - k(\omega). \qquad (3)$$

Eq. (2) is derived from Greens functions, but as the Helmholtz equation is a linear differential equation, it can be used for any solution. Compared to the 3D equation from Dean-Ben et al [22], the additional factor $\omega/(c_0 K(\omega))$ in Eq. (2) turns out to be approximately one for relevant frequencies and attenuation coefficients.

The attenuated wave in contrast to the ideal wave shows a decay in amplitude according to Eq. (2) of a factor of $\exp(-\alpha(\omega)|r|)$. As the attenuation coefficient $\alpha(\omega)$ increases with frequency, we can determine a cut-off frequency $\omega_{cut}$, for which the amplitude of the attenuated wave gets below the noise level:

$$SNR \exp(-\alpha(\omega_{cut})|r|) = 1 \; or \; \alpha(\omega_{cut}) = \frac{\ln(SNR)}{|r|}, \tag{4}$$

where the signal-to-noise-ratio $SNR$ is the amplitude of the wave without attenuation divided by the noise level and ln denotes the natural logarithm. For analyzing spatial resolution in photoacoustic imaging, the width of the acoustic signal in the time domain is essential. A small width enables high spatial resolution, which corresponds to a high frequency bandwidth. If the frequency bandwidth is limited by thermodynamic fluctuations according to Eq. (4), the spatial resolution limit according to Nyquist is half the wavelength at this frequency:

$$\delta_{resolution} = \frac{\pi}{\omega_{cut}} c(\omega_{cut}). \tag{5}$$

In time domain, according to Eq. (2), $p(r,t)$ is the inverse Fourier transform of $\tilde{p}(r,\omega)$

$$p(r,t) = \frac{1}{2\pi} \int_{-\infty}^{\infty} \tilde{p}(r,\omega) \exp(-i\omega t) \, d\omega = p_{ideal}(r,t) *_t M(r,t)$$

$$with \; M(r,t) := \frac{1}{2\pi} \int_{-\infty}^{\infty} \frac{\omega}{c_0 K(\omega)} \exp(i\gamma(\omega)|r|) \exp(-i\omega t) \, d\omega, \tag{6}$$

where $*_t$ denotes the convolution in time. The kernel M can be calculated analytically only in a few special cases, for example when $\alpha(\omega)$ is a power-law with exponent two which models liquids [10]. For discretized signals in time, Eq. (6) can be written in matrix notation:

$$\boldsymbol{p}_r = \boldsymbol{M}_r \boldsymbol{p}_{ideal}$$

$$\text{with } \boldsymbol{M}_r = \boldsymbol{F}^* \, diag\left(\frac{\omega}{c_0 K(\omega)} \exp(i\gamma(\omega)|r|)\right) \boldsymbol{F}, \tag{7}$$

where $\boldsymbol{p}_r$ and $\boldsymbol{p}_{ideal}$ are vectors and the matrix $\boldsymbol{M}_r$ describes the influence of acoustic attenuation for a propagation distance of $r$ in fat tissue. Later on, this matrix comprises also the impulse response of the piezoelectric transducer and the amplifier to enable an ideal $\delta$ - like signal for $\boldsymbol{p}_{ideal}$. Multiplication by $\boldsymbol{F}$ denotes the (discrete) Fourier transform, multiplication by its conjugate transpose $\boldsymbol{F}^*$ is the inverse Fourier transform, and $diag(\cdot)$ forms a diagonal matrix. Eq. (7) shows immediately the singular value decomposition (SVD). As discussed above the singular values decrease exponentially as $\exp(-\alpha(\omega)|r|)$. Therefore, $\boldsymbol{M}_r$ cannot be inverted exactly in a stable manner, but the pseudo-inverse matrix can be approximated using the truncated SVD (T-SVD) method (see e.g. [10]). The truncation criterion for the T-SVD method comes from the discrepancy principle and states that the inverse singular values are set to zero if they get larger than the SNR [16]. This gives the same truncation frequency as the cut-off frequency in Eq. (4), and therefore the resolution from T-SVD is the same as derived in Eq. (5). As a second reconstruction method we apply the Douglas-Rachford (DR) splitting algorithm for solving Eq. (7) with sparsity and positivity constraints.

## Results and discussion

The inverted signals using the T-SVD method are shown in Figure 2. The width of the reconstructed pulse is the reciprocal value of the cut-off-frequency of 24 MHz for 6 mm fat tissue and 11 MHz for 20 mm fat tissue. The spatial resolution according to Eq. (5) is half the width of this reconstructed pulse multiplied by the sound velocity at that frequency: 21 ns or 32 µm for 6 mm fat, and 46 ns or 70 µm for 20 mm fat, respectively.

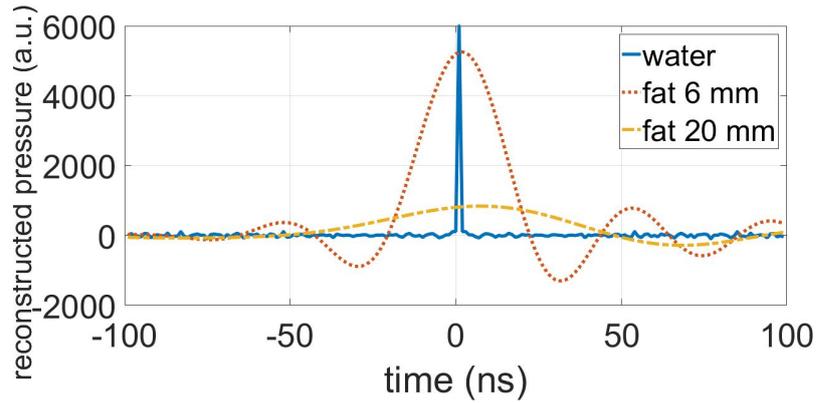

Figure 2. Reconstruction using T-SVD for regularization to compensate attenuation in fatty tissue of 6 mm thickness (red dotted line) and 20 mm thickness (yellow dashed line). This corresponds to the resolution limit from entropy production with 32 µm for 6 mm fat and 70 µm for 20 mm fat, respectively (see text). The matrix $\boldsymbol{M}_r$ was multiplied by the convolution matrix of the water signal to get a positive $\delta$ - like pulse for the pure water-signal without fatty tissue in the measurement chamber (blue solid line).

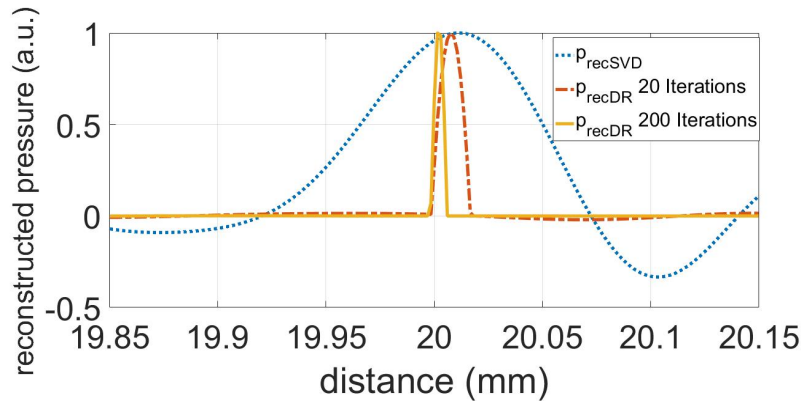

Figure 3. Compensation of acoustic attenuation of a 20 mm thick fatty tissue: T-SVD (blue dotted line), and DR method (red dashed dotted and yellow solid line) with 20 and 200 iterations, respectively. For the linear SVD reconstruction the width of the peak gives the spatial resolution, for the nonlinear DR method this is not necessarily the case. All signals are normalized to have a maximum of one.

The T-SVD method is a linear reconstruction method. If the measurement is repeated several times, it gives the same result if the measurement data is averaged and then the reconstruction is calculated, or if every measurement is reconstructed and then the average of the reconstructions is calculated. This is different for the non-linear iterative DR method. Another important difference is that for the linear method the spatial resolution of the reconstruction is given by the width of the peak, as shown in Figure 2. This is not necessarily the case for the DR method, because the width of the peak changes with the number of iterations, as shown in Figure 3 for 20 and 200 iterations. Therefore, for future work we intend to use two overlapping point-sources at a certain distance, which could be experimentally realized by flange silicon absorber having a small step in the beam-center. Due to the unfocused piezoelectric transducer, the slightly time shifted signals propagate through the fatty tissue and overlap. Resolution is then defined as the smallest distance, at which these steps can be resolved as two individual peaks. Another difficulty in defining the resolution is, that for repeated identical measurements it can happen that sometimes the peaks can be resolved and sometimes they are not resolved.

# Acknowledgement

The authors thank Gregor Thummerer from the University of Applied Sciences Upper Austria for inspiring discussions on iterative regularization methods, especially used in seismology. Measurements were funded by the project "multimodal and in-situ characterization of inhomogeneous materials" (MiCi) by the federal government of Upper Austria and the European Regional Development Fund (EFRE) in the framework of the EU-program IWB2020. Signal and data processing was also funded within the strategic economic-research program "Innovative Upper Austria 2020" of the province of Upper Austria. Parts of this work have been supported by the Austrian Science Fund (FWF), projects P 30747-N32 and P 33019-N.

# References

1   R. A. Kruger, W. L. Kiser, D. R. Reinecke, G. A. Kruger, K. D. Miller, Thermoacoustic molecular imaging of small animals, Molecular Imaging 2 (2003) 113-123.
2   P. Beard, Biomedical photoacoustic imaging, Interface Focus 1 (2011) 602-631.
3   L. V. Wang, S. Hu, Photoacoustic tomography: in vivo imaging from organelles to organs, Science 335 (2012) 1458–1462.
4   P. Burgholzer, G. J. Matt, M. Haltmeier, G. Paltauf, Exact and approximative imaging methods for photoacoustic tomography using an arbitrary detection surface, Phys. Rev. E 75 (2007) 046706.
5   P. Kuchment, L. Kunyansky, Mathematics of thermoacoustic tomography, Eur. J. Appl. Math. 19 (2008) 191–224.

6   P. Burgholzer, H. Roitner, J. Bauer-Marschallinger, H. Gruen, T. Berer, G. Paltauf, Compensation of Ultrasound Attenuation in Photoacoustic Imaging, in: Marco G. Beghi (Ed.), Acoustic Waves - From Microdevices to Helioseismology, ISBN: 978-953-307-572-3, InTech, 2011, Available from: http://www.intechopen.com/books/acoustic-waves-from-microdevices-to-helioseismology/compensation-of-ultrasound-attenuation-in-photoacoustic-imaging.
7   M. Esposito, C. Van den Broeck, Second law and Landauer principle far from equilibrium, Europhys. Lett. 95 (2011) 40004-1-6.
8   P. Burgholzer, Thermodynamic Limits of Spatial Resolution in Active Thermography, International Journal of Thermophysics 36 (2015) 2328 – 2341.
9   P. Burgholzer, G. Stockner, G. Mayr, Acoustic Reconstruction for Photothermal Imaging, Bioengineering 5 (2018) 70 1-9.
10  P. Burgholzer, J. Bauer-Marschallinger, B. Reitinger, T. Berer, Resolution Limits in Photoacoustic Imaging Caused by Acoustic Attenuation, J. Imaging 5 (2019) 13 1-11.
11  J. Bauer-Marschallinger, T. Berer, H. Gruen, H. Roitner, B. Reitinger, P. Burgholzer, Broadband high-frequency measurement of ultrasonic attenuation of tissues and liquids, IEEE Trans. Ultrason. Ferroelectr. Freq. Control. 59 (2012) 2631-2645.
12  H. Roitner, J. Bauer-Marschallinger, T. Berer, P. Burgholzer, Experimental evaluation of time domain models for ultrasound attenuation losses in photoacoustic imaging, J. Acoust. Soc. Am. 131 (2012) 3763-3774.
13  R. C. Aster, B. Borchers, C. H. Thurber, Parameter estimation and inverse problems, third ed., Elsevier, Amsterdam, 2018.
14  Y. Wang, X. Ma, H. Zhou, Y. Chen, L1-2 minimization for exact and stable seismic attenuation compensation, Geophys. J. Int. 213 (2018) 1629-1646.
15  P. Burgholzer, J. Bauer-Marschallinger, H. Grün, M. Haltmeier, G. Paltauf, Temporal back-projection algorithms for photoacoustic tomography with integrating line detectors, Inverse Problems 23 (2007) S65–S80.
16. P. Burgholzer, M. Thor, J. Gruber, G. Mayr, Three-dimensional thermographic imaging using a virtual wave concept, Journal of Applied Physics 121 (2017) 105102.
17  J. Eckstein, D. P. Bertsekas, On the Douglas—Rachford splitting method and the proximal point algorithm for maximal monotone operators, Mathematical Programming 55 (1992) 293-318.
18  J. Eckstein, W. Yao, Relative-error approximate versions of douglas– rachford splitting and special cases of the admm, Mathematical Programming 170 (2018) 417–444.
19  G. Paltauf, R. Nuster, P. Burgholzer, Weight factors for limited angle photoacoustic tomography, Physics in medicine and biology 54 (2009) 3303–3314.
20. C. B. Scruby, L. E. Drain, Laser ultrasonics techniques and applications, CRC Press, ISBN 0-7503-0050-7, 1990.
21. D. A. Hutchins, Ultrasonic generation by pulsed lasers, in: Physical Acoustics, Elsevier 18, ISBN 0893-388X, 1988, pp. 21–123.
22  X. L. Deán-Ben, D. Razansky, V. Ntziachristos, The effects of acoustic attenuation in optoacoustic signals, Phys. Med. Biol. 56 (2011) 6129–6148.
23  H. Ammari, E. Bretin, V. Jugnon, A. Wahab, Photoacoustic imaging for attenuating acoustic media, in: H. Ammari (Ed.), Mathematical Modeling in Biomedical Imaging II: Optical, Ultrasound, and Opto-Acoustic Tomographies, Lecture Notes in Mathematics 2035, Springer-Verlag, Berlin, 2012, pp. 57–84.
24  H. B. Nussenzveig, Causality and Dispersion Relations, Academic, New York, 1972.